\documentclass[12pt]{article}

\usepackage[english]{babel}

\usepackage{amsfonts}
\usepackage{amsmath}

\usepackage[hidelinks]{hyperref}

\title{Faulhaber's Problem Revisited: Alternate Methods for Deriving the Bernoulli Numbers--Part I}
\author{Christina Taylor, South Dakota School of Mines \& Technology}

\begin{document}
\maketitle
\flushleft

\subsection*{Abstract}

This paper sets the groundwork for the consideration of families of recursively defined polynomials and rational functions capable of describing the Bernoulli numbers. These families of functions arise from various recursive definitions of the Bernoulli numbers. The derivation of these recursive definitions is shown here using the original application of the Bernoulli numbers: sums of powers of positive integers, i.e., Faulhaber's problem. Since two of the three recursive forms of the Bernoulli numbers shown here are readily available in literature and simple to derive, this paper focuses on the development of the third, non-linear recursive definition arising from a derivation of Faulhaber's formula. This definition is of particular interest as it shows an example of an equivalence between linear and non-linear recursive expressions. Part II of this paper will follow up with an in-depth look at this non-linear definition and the conversion of the linear definitions into polynomials and rational functions.

\newpage
To derive expressions for the Bernoulli numbers, we return to the original context under which the Bernoulli numbers were defined: sums of powers of positive integers, or power sums, which take the form:
$$\sum_{i=1}^{n} i^m,~~~ m \in \mathbb{Z}, m \ge 0$$

Jakob Bernoulli gave the formula for these sums as:
\medskip
$$\sum_{i=1}^n i^m = \frac{1}{m+1}n^{m+1} + \frac{1}{2} n^m + \sum_{k=2}^m \frac{m!}{k!(m-k+1)!} B_k~n^{m-k+1}$$

\medskip
where $B_k$ represents the $k^{th}$ Bernoulli number. In the MAA \textit{Convergence} article ``Sums of Powers of Postive Integers," Janet Beery lists the works of several mathematicians on these sums, including both Bernoulli and Faulhaber \cite{Beery}. In her paper, Beery notes Bernoulli's derivation of these sums stemmed from the work of Pierre de Fermat, who used triangular numbers of various orders to extract formulas for the sums. Similarly, Blaise Pascal used the binomial theorem as an entry point when considering these sums. 

\medskip
In the above cases, Fermat's triangular numbers and Pascal's binomials served the same purpose: they provided a polynomial equation of arbitrary degree from which to construct formulas for power sums. Alternatively, however, Faulhaber's formula can be derived without a polynomial family by noting a recurrence relation between sums of successive powers. Using this method, Faulhaber's formula can be obtained along with, surprisingly, a definition for the Bernoulli numbers unique from that given by Bernoulli. This relation can be discovered by expanding and rearranging the terms of a sum of arbitrary degree. 

Observe the following\footnote{The resulting recursive relationship shown here is equivalent to a relationship noted by the ancient mathematician Abu Ali al-Hasan ibn al Hasan ibn al Haytham. His consideration of these sums is also detailed in Beery's article.}:
\begin{align*}
\sum_{i=1}^{n} i^{m} &= &&1^{m} + 2^{m} + 3^{m} + \cdots + n^{m}\\
&= &&1\cdot 1^{m-1} + 2\cdot 2^{m-1} + 3\cdot 3^{m-1} + \cdots + n\cdot n^{m-1}\\\\
&= &&(1^{m-1}) + (2^{m-1} + 2^{m-1}) + (3^{m-1}+ 3^{m-1} + 3^{m-1}) + \cdots + \\
& &&~~~~~~~~~~~~~~~~~~~~~~~~~~~~~~~~~~~~~~~~~~~~~~~~~~~(n^{m-1} + \cdots + n^{m-1})\\\\
&= &&( 1^{m-1} + 2^{m-1} + 3^{m-1} + \cdots + n^{m-1})~+ \\
& & &( 2^{m-1} + 3^{m-1} + \cdots + n^{m-1})~+ \\
& & & \vdots \\
& & & n^{m-1}\\
&= &&\sum_{i=1}^{n} \sum_{j=i}^{n} j^{m-1} = \sum_{i=1}^{n} \bigg( \sum_{j=1}^{n} j^{m-1} - \sum_{j=1}^{i-1} j^{m-1} \bigg)
\end{align*} 

$$\implies \sum_{i=1}^{n} i^m = \sum_{i=1}^{n} \bigg( \sum_{j=1}^{n} j^{m-1} - \sum_{j=1}^{i-1} j^{m-1} \bigg)$$

\medskip
Now suppose the function $S_k (n)$ gives the result of the sum $\sum_{i=1}^{n} i^k$. Without knowing anything of the nature of the family of the $S_k$ functions, it can be seen from the above that:

$$\sum_{i=1}^{n} i^m = \sum_{i=1}^{n} \bigg( \sum_{j=1}^{n} j^{m-1} - \sum_{j=1}^{i-1} j^{m-1} \bigg) = \sum_{i=1}^{n} \bigg( S_{m-1}(n) - S_{m-1}(i-1) \bigg)$$

$$\implies \sum_{i=1}^{n} i^m = n S_{m-1}(n) - \sum_{i=1}^{n} S_{m-1}(i-1)$$

\medskip
Knowing that the $S_k$ functions are truly polynomials, it can easily be seen from the first term in the above equation that the degree of the polynomials will increase with $m$. The second term of the above, $\sum_{i=1}^{n} S_{m-1}(i-1)$, is where the complexity of these sums and the Bernoulli numbers come in along with the inevitable introduction of the binomial theorem, which strongly flavors these sums and the resulting Bernoulli numbers.

\medskip
Ignoring Faulhaber's formula for the time being, suppose the formula for these sums was unknown and in need of discovery. From discrete mathematics, the formulas for $m = 0, 1, 2, 3$ are readily available. From these formulas, it can be hypothesized that the $S_k$ functions take the form of polynomials over the rationals of degree $k+1$. Let these polynomials be denoted using the following notation:

$$S_k(n) = \sum_{i=1}^{n} i^k = a_{k,k+1} n^{k+1} + a_{k,k} n^k + a_{k,k-1} n^{k-1} + \cdots + a_{k,1} n^{1},~~~~ a_{i,j} \in \mathbb{Q}$$

where the first subscript of the coefficients denotes the order of the sum, and the second the exponent of the term it multiplies. Note no constant term was included in these polynomials; its absence is easily justified by an empty sum with $n = 0$. While a simpler, high-level proof could be used to prove these sums are described by polynomials, an in-depth proof by strong induction can be used to acquire a definition for the coefficients.

\medskip
\medskip
\textbf{Proof of Polynomial Closed Form}

\medskip
For the base case, note by definition:

$$S_{0}(n) = \sum_{i=1}^{n} 1 = \sum_{i=1}^{n} i^0 = n$$

\medskip
For the inductive assumption, for $j = 0, 1, ..., k$ let:
$$S_j(n) = \sum_{i=1}^{n} i^j = a_{j,j+1} n^{j+1} + a_{j,j} n^j + a_{j,j-1} n^{j-1} + \cdots + a_{j,1} n^{1}, a_{i,j} \in \mathbb{Q}$$

\medskip
Using the previously stated recursive formula, the $(k+1)^{th}$ sum is given by:
$$\sum_{i=1}^{n} i^{k+1} = n S_{k}(n) - \sum_{i=1}^{n} S_{k}(i-1)$$

\newpage
For the moment, consider only the second term in the above equation. Under the inductive assumption:

$$\sum_{i=1}^{n} S_{k}(i-1) = \sum_{i=1}^{n} \big( a_{k,k+1} (i-1)^{k+1} + a_{k,k} (i-1)^{k}+ \cdots + a_{k,1} (i-1)^1\big)$$

\medskip
Expanding each of the powers of $(i-1)$ with the binomial theorem\footnote{Alternatively, the limits of summation can be changed by using a new index $j = (i-1)$; however, this route is more difficult in the end.}:

\begin{align*}
\sum_{i=1}^{n} S_{k}(i-1) &= \sum_{i=1}^{n} \bigg( a_{k,k+1} (i-1)^{k+1} + a_{k,k} (i-1)^{k}+ \cdots + a_{k,1} (i-1)^1\bigg)\\\\
&= \sum_{i=1}^{n} \bigg[~ a_{k,k+1} \bigg( \binom{k+1}{k+1}i^{k+1}(-1)^{0} + \cdots + \binom{k+1}{0}i^{0}(-1)^{k+1}\bigg)\\\\
&~~~~ + a_{k,k} \bigg( \binom{k}{k}i^{k}(-1)^{0} +\binom{k}{k-1}i^{k-1}(-1)^{1} + \cdots + \binom{k}{0}i^{0}(-1)^{k}\bigg)\\
&~~~~~ \vdots\\
&~~~~ + a_{k,2} \bigg( \binom{2}{2}i^{2}(-1)^{0} + \binom{2}{1}i^{1}(-1)^{1} + \binom{2}{0}i^{0}(-1)^{2}\bigg)\\\\
&~~~~ + a_{k,1} \bigg( \binom{1}{1}i^{1}(-1)^{0} + \binom{1}{0}i^{0}(-1)^{1}\bigg) \bigg]\\
\end{align*}

Grouping terms according to powers of $i$:
\begin{align*}
\sum_{i=1}^{n} &S_{k}(i-1) = \sum_{i=1}^{n} \bigg[~ i^{k+1}\cdot a_{k,k+1} \binom{k+1}{k+1}(-1)^{0} \\\\
&+ i^{k} \bigg( a_{k,k+1}\binom{k+1}{k}(-1)^{1} + a_{k, k}\binom{k}{k}(-1)^{0}\bigg)\\\\
&+ i^{k-1} \bigg( a_{k,k+1}\binom{k+1}{k-1}(-1)^{2} + a_{k,k}\binom{k}{k-1}(-1)^{1} + a_{k,k-1}\binom{k-1}{k-1}(-1)^{0}\bigg)\\\\[-12pt]
&~~~ \vdots \\\\[-12pt]
& + i^{1} \bigg( a_{k,k+1}\binom{k+1}{1}(-1)^{k} + \cdots + a_{k, 1}\binom{1}{1}(-1)^{0}\bigg)\\\\
& + i^{0} \bigg( a_{k,k+1}\binom{k+1}{0}(-1)^{k+1} + \cdots + a_{k, 1}\binom{1}{0}(-1)^{1}\bigg) \bigg]
\end{align*}

For ease of communication, let the coefficients for each $i^j$ with $j > 0$ be denoted as follows:

$$\alpha_{k,j} = a_{k,k+1}\binom{k+1}{j}(-1)^{k+1-j} + a_{k,k}\binom{k}{j}(-1)^{k-j} + \cdots + a_{k,j}\binom{j}{j}(-1)^{0}$$

For $j = 0$, a slightly altered formula is necessary in which the $\binom{0}{0}$ term is absent:
$$\alpha_{k,0} = a_{k,k+1}\binom{k+1}{0}(-1)^{k+1} + a_{k,k}\binom{k}{0}(-1)^{k} + \cdots + a_{k,1}\binom{1}{0}(-1)^{1}$$

\pagebreak
Returning to the derivation with these notational simplifications:

\begin{align*}
\sum_{i=1}^{n} S_{k}(i-1) &= \sum_{i=1}^{n} \big(~ \alpha_{k,k+1} \cdot i^{k+1} + \alpha_{k,k} \cdot i^{k} + \cdots \alpha_{k,1} \cdot i^{1} + \alpha_{k,0} \cdot i^{0} \big)\\
&= \alpha_{k,k+1} \sum_{i=1}^{n} i^{k+1} + \alpha_{k,k} \cdot S_{k}(n) + \cdots  + \alpha_{k,0} \cdot S_{0}(n)
\end{align*}

Substituting this back into the equation for $S_{k+1}$:

\begin{align*}
S_{k+1}(n) &= \sum_{i=1}^{n} i^{k+1} = n S_{k}(n) - \sum_{i=1}^{n} S_{k}(i-1)\\
&= n S_{k}(n) - \bigg(~\alpha_{k,k+1} \sum_{i=1}^{n} i^{k+1} + \alpha_{k,k} \cdot S_{k}(n) + \cdots  + \alpha_{k,0} \cdot S_{0}(n) ~\bigg)
\end{align*}

Solving for the $(k+1)^{th}$ sum:

$$\implies (1 + \alpha_{k,k+1})~\sum_{i=1}^{n} i^{k+1} = n S_{k}(n) - \big(~\alpha_{k,k} \cdot S_{k}(n) + \cdots  + \alpha_{k,0} \cdot S_{0}(n) ~\big)$$

$$\implies \sum_{i=1}^{n} i^{k+1} = \frac{1}{(1 + \alpha_{k,k+1})} \bigg[~n S_{k}(n) - \big(~\alpha_{k,k} \cdot S_{k}(n) + \cdots  + \alpha_{k,0} \cdot S_{0}(n) ~\big) \bigg]$$

At this point, provided $\alpha_{k,k+1} \ne -1$ to prevent division by zero, the polynomial status of the $S_k$'s has been proven. However, further manipulation is needed to obtain a definition for the polynomial coefficients, the $a_{i,j}$'s.

\medskip
Expanding each of the $S_i$'s and grouping terms according to powers of $n$:
\begin{align*}
\sum_{i=1}^{n} i^{k+1} &= \frac{1}{(1 + \alpha_{k,k+1})} \bigg[~n~\big(a_{k,k+1} \cdot n^{k+1} + a_{k,k} \cdot n^{k} + \cdots + a_{k,1} \cdot n^{1}\big) \\\\
&~~~~~~~~ - \alpha_{k,k} \big(a_{k,k+1} \cdot n^{k+1} + a_{k,k} \cdot n^{k} + \cdots a_{k,1} \cdot n^{1}\big)\\\\
&~~~~~~~~ - \alpha_{k,k-1} \big(a_{k-1,k} \cdot n^{k} + a_{k-1,k} \cdot n^{k-1} + \cdots a_{k-1,1} \cdot n^{1}\big)\\\\[-12pt]
&~~~~~~~~ \vdots \\\\[-12pt]
&~~~~~~~~ - \alpha_{k,1} \big( a_{1,2} n^2 + a_{1,1} n^1\big) - \alpha_{k,0} \big( a_{0,1} n^1 \big) ~\bigg]
\end{align*}

\begin{align*}
\sum_{i=1}^{n} i^{k+1} = \frac{1}{(1 + \alpha_{k,k+1})} \bigg[&~ n^{k+2} \cdot a_{k,k+1} ~+ \\
&n^{k+1} ~\big(a_{k,k} - \alpha_{k,k} \cdot a_{k,k+1}\big) ~+ \\
&n^{k} ~\big(a_{k,k-1} - \alpha_{k,k} \cdot a_{k,k} - \alpha_{k,k-1}\cdot a_{k-1,k}~\big) ~+ \\
& \vdots \\
&n^{2} ~\big(a_{k,1} - \alpha_{k,k} \cdot a_{k,2} - \alpha_{k,k-1} \cdot a_{k-1,2} - \cdots -\alpha_{k,1} \cdot a_{1,2} \big) ~- \\
&n^{1} ~\big(~\alpha_{k,k} \cdot a_{k,1} + \alpha_{k,k-1} \cdot a_{k-1,1} + \cdots + \alpha_{k,0} \cdot a_{0,1}~\big)~\bigg]
\end{align*}

With the above, the proof is now complete, and a generalized definition for the coefficients of the polynomials can be observed to be:
\begin{multline*}
a_{m,k} = \frac{1}{1 + \alpha_{m-1,m}}\bigg[~a_{m-1,k-1} - \big( \alpha_{m-1,m-1} \cdot a_{m-1,k} ~+  \\
\alpha_{m-1,m-2} \cdot a_{m-2,k} + \cdots + \alpha_{m-1,k-1} \cdot a_{k-1,k} \big)~\bigg]
\end{multline*}

This indexing, however, is not actually the most natural setup. Given that the coefficients' recursive definitions become stronger the lower the order of the term, it is more convenient to index the coefficients by offsetting from the leading term, or rather its consort, $n^m$. Put more mathematically\footnote{Note this offset-style indexing corresponds to the alternative indexing starting at $-1$ sometimes used when dealing with the Bernoulli numbers.}:
\begin{multline*}
a_{m,m-x} = \frac{1}{1 + \alpha_{m-1,m}}\bigg[~a_{m-1,m-x-1} - \big( \alpha_{m-1,m-1} \cdot a_{m-1,m-x} ~+  \\
\alpha_{m-1,m-2} \cdot a_{m-2,m-x} + \cdots + \alpha_{m-1,(m-1)-x} \cdot a_{(m-1) -x,m-x} \big)~\bigg]
\end{multline*}

Using this indexing, the above definition is valid for $ 0 \le x \le m-2$, i.e., all but the coefficients on the first ($x = -1, k = m+1$) and last ($x = m - 1, k = 1$) terms, whose definitions are given by:

$$a_{m,m+1} = \frac{a_{m-1,m}}{1+\alpha_{m-1,m}}$$
$$a_{m,1} = \frac{-1}{1 + \alpha_{m-1,m}}\bigg[~\alpha_{m-1,m-1} \cdot a_{m-1,1} + \alpha_{m-1,m-2} \cdot a_{m-2,1} + \cdots + \alpha_{m-1,0} \cdot a_{0,1}~\bigg]$$

Using the above, closed forms for the coefficients can be obtained. For the sake of some sense of brevity, the process used to do so will be described below but not shown in detail.

\medskip
To obtain closed forms for the coefficients, note that the coefficients are multivariate: they are functions of both $x$ and $m$. Additionally, they are strongly recursive. However, holding $x$ constant creates sequences solely in terms of $m$. The base cases for each of these offset-sequences are attainable thanks to the special definitions for $\alpha_{m,0}$ and $a_{m,1}$. Finding a closed form for the leading term, $x = -1$, is a trivial exercise as it depends only on the previous term in the $x = -1$ sequence\footnote{This sequence also puts to rest the previous concern for division by zero previously noted as its terms are always positive.}; for all other offsets, the above strong recursive definitions must first be weakened, i.e., its dependence on sequences with offsets different from its own pruned.

\medskip
Suppose the sequence to be weakened was $x = b$. To weaken the recursive definition, the closed form for sequences with $x < b$ must be known; these closed forms are then to be substituted into the strong recursive definitions given above to obtain an expression solely in terms of $m, x,$ and previous terms of the $x = b$ sequence. 

\pagebreak
It can be shown that these weakened/condensed recursive forms in general take the form:

$$a_{m,m-x} = a_{m-1,(m-1)-x}\frac{m~(m-x-1)}{(m+1)(m-x)} + \frac{C_x \cdot m!}{(m+1) (m-x)!}$$

\medskip
with base cases given by the independent term:
$$a_{x+1,1} = \frac{C_x \cdot (x+1)!}{x+2}$$

where the $C_{i}$'s are constants unique to each offset sequence (these eventually lead to the Bernoulli numbers). This condensed recursive form simplifies the proof of the closed form of the coefficients, which at this stage are given by:

$$a_{m,m-x} = \frac{C_x \cdot m!}{(x+2)(m-x)!}$$

\medskip
To conduct this proof, proof by mathematical induction is first used to show the condensed recursive definition yields the closed form above. Secondly, a proof by strong induction then proves that the above closed form, when substituted into the strongly recursive definition, yields the condensed recursive definition. These two proofs/sub-proofs also yield the definition for the $C_i$'s:

$$C_x = - \bigg[~\frac{C_{x-1}}{x+1}\beta_{0} + \frac{C_{x-2}}{x}\beta_{1} + \cdots + \frac{C_{0}}{2}\beta_{x-1} + \frac{C_{-1}}{1}\beta_{x}^* ~\bigg],~~~~ C_{-1} = 1$$

The $\beta_{i}$ terms arise from the $\alpha_{k,i}$'s and are given by:

$$\beta_{x} = \frac{C_{-1}(-1)^{x+1}}{1 \cdot (x+1)!} + \frac{C_{0}(-1)^{x}}{2 \cdot x!} + \cdots + \frac{C_{x-1}(-1)^{1}}{(x+1) \cdot 1!} + \frac{C_{x}(-1)^{0}}{(x+2)\cdot 0!}$$

$$\beta_{x}^* = \frac{C_{-1}(-1)^{x+1}}{1 \cdot (x+1)!} + \frac{C_{0}(-1)^{x}}{2 \cdot x!} + \cdots + \frac{C_{x-2}(-1)^{2}}{x \cdot 2!} + \frac{C_{x-1}(-1)^{1}}{(x+1)\cdot 1!}$$

\pagebreak
The $C_i$'s are a few steps away from the Bernoulli numbers. It can be observed and proven that a given $C_x$ has a factor of one over $(x+1)!$, which leads to the $D_x$ sequence:

$$C_x = \frac{D_x}{(x+1)!} \implies D_x = (x+1)! C_x$$

At this point, using the base case of $\sum_{i=1}^{n} i^{0} = n$ and proof by strong and mathematical induction, this derivation has proven that:

$$\sum_{i=1}^{n} i^{m} = \frac{D_{-1}~ m!}{1!~(m+1)!}n^{m+1} + \frac{D_{0}~ m!}{2!~m!}n^{m} + \cdots + \frac{D_{m-1}~ m!}{(m+1)!~1!}n^{1}$$
$$\implies \sum_{i=1}^{n} i^{m} = \sum_{x = -1}^{m-1} \frac{D_{x}~ m!}{(x+2)!~(m-x)!}n^{m-x}$$

\medskip
which is completely analogous to Faulhaber's formula. From the above it can be observed that when the traditional indexing from 0 is used for the Bernoulli numbers, they are related quite simply to the $D_x$ sequence by:

$$D_x = \frac{B_{x+1}}{x+2} \implies B_{x} = (x+1) D_{x-1}$$

\medskip
Alternatively, when the Bernoulli numbers are indexed from $-1$:

$$D_x = \frac{B_x}{x+2} \implies B_x = (x+2) D_x$$

The most interesting consequence of this equivalency, however, is the difference in the definitions of the Bernoulli numbers and the $D_x$'s. There are many definitions for the Bernoulli numbers, but the typical recursive form is linear with respect to the $B_i$'s; the definition for the $D_x$'s, however, is non-linear with respect to the other $D_i$'s. Ignoring the Bernoulli numbers and their many definitions, this duplicity in definition can be observed using only the $D_x$'s. Using what has been proven thus far, consider the case of a power sum with $n=1$ for any $m$:

$$\sum_{i=1}^1 i^m = 1^m = 1 = \frac{D_{-1}~ m!}{1!~(m+1)!}(1)^{m+1} + \frac{D_{0}~ m!}{2!~m!}(1)^{m} + \cdots +  \frac{D_{m-1}~ m!}{(m+1)!~1!}(1)^{1}$$

$$ \implies 1 = \frac{D_{-1}~ m!}{1!~(m+1)!} + \frac{D_{0}~ m!}{2!~m!} + \cdots + \frac{D_{m-2}~ m!}{m!~2!} + \frac{D_{m-1}~ m!}{(m+1)!~1!}$$

\medskip
For the sake of notation, let $m$ become $(x+1)$ and then solve for $D_x$:

$$D_x = 1 - (x+1)!\bigg[~ \frac{D_{-1}}{1!~(x+2)!} + \frac{D_{0}}{2!~(x+1)!} + \cdots + \frac{D_{x-2}}{x!~3!} + \frac{D_{x-1}}{(x+1)!~2!} ~\bigg]$$

\medskip
This definition for the $D_x$ is analogous to one of the recursive definitions commonly given for the Bernoulli numbers, but recall this derivation gives a non-linear definition for the $D_x$'s:

$$D_x = - \bigg[~ \frac{D_{x-1}}{(x+1)!} \beta_{0} + \frac{D_{x-2}}{x!} \beta_{1} + \cdots + \frac{D_{0}}{2!} \beta_{x-1} + \frac{D_{-1}}{1!} \beta_{x}^*~\bigg]$$

\medskip
$$\beta_{x} = \frac{D_{-1}(-1)^{x+1}}{1! \cdot (x+1)!} + \frac{D_{0}(-1)^{x}}{2! \cdot x!} + \cdots + \frac{D_{x-1}(-1)^{1}}{(x+1)! \cdot 1!} + \frac{D_{x}(-1)^{0}}{(x+2)!\cdot 0!}$$

\medskip
$$\beta_{x}^* = \frac{D_{-1}(-1)^{x+1}}{1! \cdot (x+1)!} + \frac{D_{0}(-1)^{x}}{2! \cdot x!} + \cdots + \frac{D_{x-2}(-1)^{2}}{x! \cdot 2!} + \frac{D_{x-1}(-1)^{1}}{(x+1)!\cdot 1!}$$

\medskip
Additionally, there is another linear definition for the $D_x$'s that can be gleaned from this derivation using the idea of an empty sum where $n = -1$:

$$\sum_{i=1}^{-1} i^m = 0 = \frac{D_{-1}~ m!}{1!~(m+1)!}(-1)^{m+1} + \frac{D_{0}~ m!}{2!~m!}(-1)^{m} + \cdots + \frac{D_{m-1}~ m!}{(m+1)!~1!}(-1)^{1}$$

$$D_x = (x+1)!\bigg[~ \frac{D_{-1}(-1)^{x+2}}{1!~(x+2)!} + \frac{D_{0}(-1)^{x+1}}{2!~(x+1)!} + \cdots + \frac{D_{x-2}(-1)^{3}}{x!~3!} + \frac{D_{x-1}(-1)^{2}}{(x+1)!~2!} ~\bigg]$$

\medskip
These contending definitions (among others) testify of the fascinating nature of the Bernoulli numbers. From the above, a family of recursive polynomials and rational functions can be defined as can a matrix equation. The matrix equation and its consequences have been detailed by Giorgio Pietrocola and thus will not be discussed in great depth in Part II of this paper \cite{Pietrocola}. The conversion of the linear recursive forms into a family of polynomials, dubbed the partial polynomials, and a similarly defined family of rational functions, however, will be discussed in detail in Part II of this paper, as will the equivalence of the non-linear and linear definitions.


\newpage
\begin{thebibliography}{10}
\bibitem{Beery} Beery, Janet. ``Sums of Powers of Positive Integers." \textit{Convergence}. MAA. July 2010.

\bibitem{Pietrocola} Pietrocola, Giorgio. ``On polynomials for the calculation of sums of powers of successive integers and Bernoulli numbers..." June 2017. \scriptsize{\url{http://www.pietrocola.eu/EN/Theoremsonthesumofpowersofsuccessiveintegersbygiorgiopietrocola%20.pdf}}
\end{thebibliography}
\end{document}